\documentclass{amsart}
\usepackage{amscd} 
\usepackage{amsfonts} 
\usepackage{amssymb} 
\usepackage{latexsym} 
\input{diagrams}
\newcommand{\ncm}{\newcommand}

\newtheorem{theorem}{Theorem}[section]
\newtheorem{prop}[theorem]{Proposition}
\newtheorem{lemma}[theorem]{Lemma}
\newtheorem{cor}[theorem]{Corollary}
\newtheorem{lem&def}[theorem]{Lemma \& Definition}
\newtheorem{definition}[theorem]{Definition}

\ncm{\End}{\mbox{\rm End}\,}

\def\Hom{\mbox{\rm Hom}\,}

\def\id{\mbox{\rm id}}

\def\into{\hookrightarrow}
\def\to{\rightarrow}

\def\o{\otimes}    
\def\bra{\langle}
\def\ket{\rangle}

\ncm{\rarr}[1]{\stackrel{#1}{\longrightarrow}}
\ncm{\larr}[1]{\stackrel{#1}{\longleftarrow}}
\def\cop{\Delta}

\def\eps{\varepsilon}

\def\du1{\hat 1}

\def\-1{_{(-1)}}
\def\0{_{(0)}}
\def\1{_{(1)}}
\def\2{_{(2)}}
\def\3{_{(3)}}

\def\|{\, | \,}

\def\du1{\hat 1}

\begin{document}

\title{A note on Galois theory for bialgebroids}
\author{Lars Kadison}
\address{Matematiska Institutionen \\ G{\" o}teborg 
University \\ 
S-412 96 G{\" o}teborg, Sweden} 
\email{lkadison@c2i.net} 
\date{}
\thanks{The author thanks A.A.~Stolin and K.~Szlach\'anyi 
for discussions.}
\subjclass{16W30 (13B05, 20L05, 16S40, 81R50)}  
\date{} 

\begin{abstract} 
 In this note we reduce certain proofs in \cite{KS, Karl, AMA} to depth two quasibases from one side only.  This minimalistic approach leads to a characterization of Galois extensions for finite
projective bialgebroids without the Frobenius extension property:
a proper algebra extension is a left $T$-Galois extension for some right finite projective
left bialgebroid $T$ over some algebra $R$ if and only if it is of left depth two and left balanced. Exchanging left and right in this statement,
we have also a characterization of right Galois extensions for
left finite projective right bialgebroids. As a corollary,
we obtain insights into split monic Galois mappings and
endomorphism ring theorems for depth two extensions.   
\end{abstract} 
\maketitle

\section{Introduction and Preliminaries}

Hopf algebroids arise as the endomorphisms of fiber functors from certain tensor categories
to a bimodule category over base algebra. For example, 
Hopf algebroids over a one-dimensional base algebra are Hopf algebras
while Hopf algebroids over a separable $K$-algebra base are  weak Hopf algebras.  
Galois theory for right or left bialgebroids were introduced recently in \cite{Karl, AMA, nov} based on  
the theory of Galois corings \cite{BW} and ordinary definitions of Galois extensions \cite{Mo, CDG} 
with applications to depth two extensions.  In \cite{AMA} Frobenius extensions that are right Galois over
a left finite projective right bialgebroid are characterized as being of depth two and right balanced.  
Then a Galois theory for Hopf algebroids, especially of Frobenius type, was introduced
in \cite{Bo, BaS} with applications to Frobenius extensions of depth two and weak Hopf-Galois extensions over
finite dimensional quantum groupoids (among other things, showing that these are Frobenius extensions).
Although they break with the tradition of defining Galois extensions over bialgebras
and have a more complex definition, Galois extensions over Hopf algebroids have more properties in common 
with Hopf-Galois extensions.  However, several of these properties will follow from any Galois theory for bialgebroids
which is in possession of two Galois mappings equivalent due to a bijective antipode, sometimes denoted by $\beta$ and $\beta'$, as is the case for finite Hopf-Galois theory \cite[ch.\ 8]{Mo}, finite weak Hopf-Galois theory \cite[section 5]{nov},
or possibly some future, useful weakening of Hopf-Galois theory to Hopf algebroids over a symmetric algebra,
Frobenius algebra or 
some other type of base algebra.  

In \cite{BaS} a characterization similar to \cite{AMA} of depth two Frobenius one-sided balanced extensions are
given in terms of Galois extensions over Hopf algebroids with integrals. This shows in a way that
the main theorem in \cite{AMA} makes no essential use of the hypothesis of Frobenius extension
(only that a Frobenius extension is of left depth two if and only if it is of right depth two),
and it would be desirable at some time to remove the Frobenius hypothesis.  This is then the objective of this
paper:  to show that Galois extensions over one-sided finite projective bialgebroids are characterized by
one-sided depth two and balance conditions on the extension (Theorem~\ref{th-characterization}). 
This requires among other things some care in re-doing the two-sided arguments in \cite{KS} to show that the structure $ T := (A \o_B A)^B$
on a one-sided depth two extension $A \| B$ with centralizer $R$ is still a one-sided finite projective
right bialgebroid (proposition~\ref{prop-tee}). In an appendix, we do a similar one-sided derivation
of the left bialgebroid structure on the $R$-dual $S := \End {}_BA_B$. These two sections may be read as
an introduction to depth two theory.

Let $K$ be any commutative ground ring in this paper. All algebras
are unital associative $K$-algebras and modules over these are symmetric unital $K$-modules. 
We say that $A \| B$ is an \textit{extension} (of algebras) if there is an algebra homomorphism $B \to A$, \textit{proper} if this
is monic.  This homomorphism induces the natural bimodule
structure ${}_BA_B$ which is most important to our set-up.
The extension $A \| B$ is \textit{left depth two} (left D2) if 
the tensor-square $A \o_B A$ is centrally projective w.r.t.\
$A$ as natural $B$-$A$-bimodules: i.e., 
$$ {}_BA \o_B A_A \oplus * \cong \oplus^n {}_BA_A. $$
This last statement postulates the existence then of a
split $B$-$A$-epimorphism from $A \o_B A$ to direct product
of $A$ with itself $n$ times.  

Making the clear-cut identifications  $\Hom ({}_BA \o_B A_A, {}_BA_A)
\cong $ $\End {}_BA_B$ and $\Hom ({}_BA_B, {}_BA \o_B A_A) \cong $
$(A \o_B A)^B$, we see that left D2 is characterized by
there being  left D2 quasibases $t_i \in (A \o_B A)^B$
and $\beta_i \in \End {}_BA_B$ such that for all $a, a' \in A$
\begin{equation}
\label{eq: left D2 quasibase}
a \o_B a' = \sum_{i=1}^n t_i \beta_i(a)a'.
\end{equation}
The algebras $\End {}_BA_B$ and $(A \o_B A)^B$ (note
that the latter is isomorphic to $\End {}_A A \o_B A_A$
and thus receives an algebra structure) are so important
in depth two theory that we fix (though not unbendingly) brief notations for these:
$$ S := \End {}_BA_B \ \ \ \ \ \ \  T := (A \o_B A)^B. $$

Similarly, a right depth two extension $A \| B$ is defined by switching from the natural $B$-$A$-bimodules in the definition
above to the natural
$A$-$B$-bimodules on the same structures.  Thus an extension
$A \| B$ is \textit{right D2} if ${}_A A \o_B A_B \oplus * \cong \oplus^m {}_AA_B$.  Equivalently, if there are
$m$ paired elements $u_j \in T$, $\gamma_j \in S$ such that
\begin{equation}
\label{eq: right D2 quasibase}
a \o a' = \sum_{j=1}^m a \gamma_j(a') u_j
\end{equation}
for all $a, a' \in A$.  

A depth two extension is one that is both left and right D2.
These have been studied in \cite{KS,Karl,AMA} among others, 
but without a focus on left or right D2 extensions.

Let $t, t'$ be elements in $T$, where we write $t$
in terms of its components using a notation that
suppresses a possible summation in $A \o_B A$:
$t = t^1 \o t^2$. 
Then the algebra structure on $T$ is simply
\begin{equation}
tt' = {t'}^1 t^1 \o t^2 {t'}^2, \ \ \ \ \ 1_T = 1_A \o 1_A
\end{equation}

There is a standard ``groupoid'' way to produce right and left
bialgebroids, which we proceed to do for $T$.  
There are two commuting embeddings of $R$ and its opposite algebra
in $T$.  A ``source'' mapping $s_R: R \to T$ given by
$s_R(r) = 1_A \o r$, which is an algebra homomorphism.
And a ``target'' mapping $t_R: R \to T$ given by
$t_R(r) = r \o 1_A$ which is an algebra anti-homomorpism
and clearly commutes with the image of $s_R$.  Thus it makes
sense to give $T$ an $R$-$R$-bimodule structure via $s_R$, $t_R$
from the right:  $r \cdot t \cdot r' = $ $t s_R(r') t_R(r) = $
$t(r \o r') = $ $rt^1 \o t^2 r'$, i.e., ${}_RT_R$ is given by
\begin{equation}
\label{eq: bimod}
r\cdot t^1 \o t^2 \cdot r' = rt^1 \o t^2 r'
\end{equation}

\begin{prop}
\label{prop-tee}
Suppose $A \| B$ is either a right or a left D2 extension.
Then $T$ is a right $R$-bialgebroid, which is either
left f.g.\ $R$-projective or right f.g.\ $R$-projective
respectively.
\end{prop}
\begin{proof}
First we suppose $A \| B$ is left D2 with quasibases
$t_i \in T$, $\beta_i \in S$.  The proof that $T$ is
a right $R$-bialgebroid in \cite[5.1]{KS} carries through
verbatim except in one place where  right D2 quasibases made
a brief appearance, where coassociativity of the coproduct
needs to be established through the introduction of an isomorphism.  Thus we need to see that
$$ T \o_R T \o_R T \stackrel{\cong}{\longrightarrow} 
(A \o_B A \o_B A \o_B A)^B$$
via $t \o t' \o t'' \mapsto t^1 \o t^2 {t'}^1 \o {t'}^2 {t''}^1 \o {t''}^2$.  
The  inverse is given by
$$ a_1 \o a_2 \o a_3 \o a_4 \longmapsto 
\sum_{i,j} t_i \o_R t_j \o_R (\beta_j(\beta_i(a_1)a_2)a_3 \o_B a_4). $$
for all $a_i \in A$ ($i = 1,2,3,4$).

In the case that we only use right D2 quasibases,
this inverse is given by
\begin{equation}
\label{eq: coassoc}
a_1 \o a_2 \o a_3 \o a_4 \longmapsto
\sum_{j,k} a_1 \o a_2 \gamma_k(a_3 \gamma_j(a_4)) \o_R u_k \o_R u_j.
\end{equation}
Both claimed inverses are easily verified as such by using the
right and left D2 quasibase equations repeatedly.

The module $T_R$ is finite projective since eq.~(\ref{eq: left D2
quasibase}) implies a dual bases equation $t =$ $ \sum_i t_i f_i(t)$, 
for each $t \in T \subseteq A\o_B A$, 
 where $f_i(t) := \beta_i(t^1)t^2$ define $n$ maps in
$\Hom (T_R, R_R)$.

Suppose $A \| B$ is right D2 with quasibases $u_j \in T$, $\gamma_j \in
S$.  The algebra structure on $T$ is given in the introduction
above as is the $R$-$R$-bimodule structure.
What remains is specifying the $R$-coring structure on $T$
and checking the five axioms of a right bialgebroid.
The coproduct $\cop: T \to T \o_R T$ is given
by 
\begin{equation}
\label{eq: cop}
\cop(t) := \sum_j (t^1 \o_B \gamma_j(t^2)) \o_R u_j,
\end{equation}
 which is clearly left $R$-linear, and right $R$-linear as well 
since
$$ \cop(tr) = \sum_j t^1 \o \gamma_j(t^2 r) \o u_j \stackrel{\cong}{\longmapsto} t^1 \o 1 \o t^2 r $$
under the isomorphism $T \o_R T \cong (A \o_B A \o_B A)^B$
given by $t \o t' \mapsto $ $t^1 \o t^2 {t'}^1 \o {t'}^2$,
which is identical to the image of 
$$ \cop(t)r = \sum_j t^1 \o \gamma_j(t^2) \o u_j r \mapsto 
t^1 \o 1 \o t^2 r .$$

Coassociativity $(\cop \o \id_T) \cop = (\id_T \o \cop) \cop$ follows
from applying the isomorphism $$T \o_R T \o_R T \cong
(A \o_B A \o_B A \o_B A)^B$$ given above in this proof 
to the left-hand and right-hand sides applied to a $t \in T$: 
$$ \sum_j \cop(t^1 \o \gamma_j(t^2)) \o u_j = \sum_{j,k}
t^1 \o \gamma_k(\gamma_j(t^2)) \o_R u_k \o_R u_j 
\stackrel{\cong}{\longmapsto} t^1 \o_B 1_A \o_B 1_A \o_B t^2. $$
$$ \sum_j (t^1 \o \gamma_j(t^2)) \o_R \cop(u_j)
= \sum_{j,k} (t^1 \o \gamma_j(t^2)) \o_R (u_j^1 \o \gamma_k(u_j^2))
\o_R u_k $$
which also maps into $t^1 \o_B 1_A \o_B 1_A \o_B t^2$ under the
same isomorphism.

The counit $\eps: T \to R$ of the $R$-coring $T$ is given by
\begin{equation}
\eps(t) := t^1 t^2
\end{equation}
i.e., the multiplication mapping $A \o_B A \to A$ restricted
to $T$ (and taking values in $R$ since $bt = tb$ for all $b \in B$).
Clearly, $\eps(rtr') = r \eps(t) r'$ for $r,r' \in R$, $t \in T$,
and that $(\id_T \o_R \eps)\cop = \id_T$ $= (\eps \o_R \id_T)\cop$
since $\sum_j t^1 \gamma_j(t^2)u_j = t$
and $\sum_j \gamma_j(a)u_j^1u_j^2 = a$ for $a \in A$. 
This shows $(T,R,\cop, \eps)$ is an $R$-coring.

We next verify the five axioms of a right bialgebroid \cite[2.1]{KS}.

\begin{enumerate}
\item $\cop(1_T) = 1_T \o 1_T$ since $\gamma_j(1_A) \in R$
and $1_T = \sum_j \gamma_j(1_A)u_j$.

\item $\eps(1_T) = 1_A$ since $1_T = 1_A \o 1_A$.

\item $\eps(tt') = \eps(t_R(\eps(t)) t') = \eps(s_R(\eps(t)) t')$
($t,t' \in T$) since $\eps(tt') =$ $ {t'}^1 t^1 t^2 {t'}^2$, 
$t_R(\eps(t)) =$ $ t^1 t^2 \o_B 1_A$ and $s_R(\eps(t)) =$ $ 1_A \o_B
t^1 t^2 $.

\item $(s_R(r) \o 1_T)\cop(t) = (1_T \o t_R(r))\cop(t)$ for
all $r \in R$, $t \in T$ since the left-hand side is
$$\sum_j (t^1 \o_B r \gamma_j(t^2)) \o_R u_j \stackrel{\cong}{\longmapsto} t^1 \o_B r \o_B t^2 $$
under the isomorphism $T \o_R T \cong (A \o_B A \o_B A)^B$
given by $t \o_R t' \mapsto$ $ t^1 \o_B t^2 {t'}^1 \o_B {t'}^2$
and the right-hand side is equal to
$$ \sum_j (t^1 \o_B \gamma_j(t^2)) \o_R (u^1_j r \o_B u^2_j) 
\stackrel{\cong}{\longmapsto} t^1 \o_B r \o_B t^2 $$ 
with the same image element.  

\item $\cop(tt') = \cop(t)\cop(t')$ for all $t,t' \in T$ in the tensor subalgebra (denoted by $T \times_R T$ with the straightforward
tensor multiplication)
of $T \o_R T$ (which makes sense thanks to the previous axiom). 
This follows from both sides having the image element ${t'}^1 t^1 \o 1_A
\o t^2 {t'}^2 $ under the isomorphism $T \o_R T \cong $
$(A \o_B A \o_B A)^B$, which is clear for the left-hand side of the equation
and for  the right-hand side we note it equals 
$$ \sum_{j,k} ({t'}^1 t^1 \o_B \gamma_j(t^2) \gamma_k({t'}^2) )
\o_R (u_k^1 u_j^1 \o_B u_j^2 u_k^2 ). $$
Now apply $t \o t' \mapsto t^1 \o_B t^2 {t'}^1 \o_B {t'}^2$
and the right D2 quasibase equation twice. 
\end{enumerate}
This completes the proof that $(T, R, s_R, t_R, \cop, \eps)$
is a right bialgebroid. 

Finally ${}_RT$ is finite projective via an application
of the right D2 quasibase eq.~(\ref{eq: right D2 quasibase}).
\end{proof}

A right comodule algebra is an algebra in the tensor category of right $R$-comodules
\cite{BaS}.  In detail, the definition is equivalent to the following.   
  
\begin{definition}
Let $T$ be any right bialgebroid $(T, R, \tilde{s}, \tilde{t},$ $ \cop, $ $ \eps)$ over any base algebra $R$.   
A right $T$-comodule algebra $A$ is an algebra $A$ with algebra homomorphism $R \to A$ (providing the $R$-$R$-bimodule structure on $A$)
 together with a coaction $\delta: A \to A \o_{R} T$, where  values $\delta(a)$ are denoted by the Sweedler
notation $a\0 \o a\1$, 
such that $A$ is a right $T$-comodule over the $R$-coring $T$ \cite[18.1]{BW}, 
$ \delta(1_{A}) = 1_{A} \o 1_{T} $, $ra\0 \o a\1 =$ $ a\0 \o \tilde{t}(r)a\1$
for all $r \in R$,
and $\delta(aa') = \delta(a) \delta(a')$ for all $a,a' \in A$.   The subalgebra of coinvariants
is ${A}^{\rm co \, T} := \{ a \in A | \delta(a) = a \o 1_{T} \}$.
We also call $A$ a right $T$-extension of ${A}^{\rm co \, T}$.   
\end{definition} 

\begin{lemma}
For the right $T$-comodule $A$ introduced just above,
 $R$ and ${A}^{\rm co \, T}$ commute in $A$.
\end{lemma}
\begin{proof}
We note that 
$$\rho(rb) = b \o_R \tilde{s}(r) = \rho(br) $$
for $r \in R$, $b \in {A}^{\rm co \, T}$. But $\rho$
is injective by the counitality of comodules, so
$rb = br$ in $A$ (suppressing the morphism $R \to A$).  
\end{proof}

\begin{definition}
Let $T$ be any  right bialgebroid
over any algebra $R$.  
A  $T$-comodule algebra $A$ is a right $T$-Galois extension of its coinvariants
$B$ if the (Galois) mapping 
$\beta: A \o_{B} A \to A \o_{R} T$ defined
by $\beta(a \o a') = a{a'}\0 \o {a'}\1 $ is bijective.  
\end{definition} 

Left comodule algebras over left bialgebroids and their left Galois extensions
are defined similarly, the details of which are in \cite{nov}.  The values of the coaction is in this case denoted
by $a\-1 \o a\0$ and the Galois mapping by $a \o a' \mapsto a\-1 \o a\0 a'$.

\section{A characterization of Galois extensions for bialgebroids}

We recall that a module ${}_AM$ is \textit{balanced} if all 
the endomorphisms of the natural module $M_E$ where $E = \End {}_AM$
are uniquely left multiplications by elements of $A$:
$A \stackrel{\cong}{\longrightarrow} \End M_E$ via $\lambda$.  
In particular, ${}_AM$ must be faithful.
 
\begin{theorem}
\label{th-characterization}
Let $A \| B$ be a proper algebra extension. Then
\begin{enumerate}
\item  $A \| B$ is a right $T$-Galois extension for
some left finite projective right bialgebroid $T$ over
some algebra $R$ if and only if $A \| B$ is right D2 and right balanced.
\item $A \| B$ is a left $T$-Galois extension for
some right finite projective left bialgebroid $T$ over
some algebra $R$ if and only if $A \| B$ is left D2 and left balanced.
\end{enumerate}
\end{theorem}
\begin{proof}
($\Rightarrow$)
Suppose $T$ is a left finite projective right bialgebroid
over some algebra $R$. 
Since  ${}_RT \oplus * \cong {}_RR^t$ for some positive integer $t$,  we apply to this 
the functor $A \o_R -$ from left $R$-modules into $A$-$B$-bimodules which results in 
${}_AA\! \o_B\! A_B \oplus * \cong {}_AA_B^t$, after using the Galois $A$-$B$-isomorphism 
$A \o_B A \cong  $ $A \o_R T$. Hence, $A | B$ is right D2.

Let $\mathcal{E} := \End A_B$.  We show $A_B$ is balanced by
the following device. 
Let $R$ be an algebra, $M_R$ and ${}_RV$ modules with ${}_RV$ finite projective.  If 
$\sum_j m_j \phi(v_j) = 0$ for all $\phi$ in the left $R$-dual ${}^*V$, then $\sum_j m_j \o_R v_j = 0$.
This follows immediately by using dual bases $f_i \in {}^*V$, $w_i \in V$.  

Given $F \in \End {}_{\mathcal{E}}A$, it suffices to show that $F = \rho_b$ for some
$b \in B$.  Since $\lambda_a \in \mathcal{E}$, $F \circ \lambda_a = \lambda_a \circ F$ for all $a \in A$,
whence $F = \rho_{F(1)}$.  Designate $F(1) = x$.  If we show that $x\0 \o x\1 = x \o 1$ after applying
the right $T$-valued coaction on $A$,
then $x \in A^{\rm co \, T} = B$. For each $\alpha \in \Hom ({}_RT, {}_RR)$, define $\overline{\alpha} \in \End A_B$
by $\overline{\alpha} (y) = y\0 \alpha(y\1)$. Since $\rho_r \in \mathcal{E}$ for each $r \in R$ by lemma, 
$$x \alpha(1_T) = F(\overline{\alpha} (1_A)) = \overline{\alpha} (F(1_A)) = x\0 \alpha(x\1) $$
for all $\alpha \in {}^*T$. Hence $x\0 \o_R x\1 = x \o 1_T$.

($\Leftarrow$) It follows from the proposition
that a right D2 extension $A \| B$ has a left finite projective
right bialgebroid $T := (A \o_B A)^B$ over the centralizer $R$
of the extension. Let $R \into A$ be the inclusion mapping.
We check that $A$ is a right $T$-comodule algebra via the
coaction $\rho_R: A
\to A \o_R T$ on $A$ given by
\begin{equation}
\label{eq: coaction}
\rho_R(a) = a\0 \o a\1 := \sum_j \gamma_j(a) \o u_j.
\end{equation}
First, we demonstrate several properties by using the isomorphism
$\beta^{-1}: A \o_R T \stackrel{\cong}{\longmapsto} A \o_B A$
given by $\beta^{-1}(a \o t) = at = at^1 \o t^2$ \cite[3.12(iii)]{KS}
with inverse $\beta(a \o a')= \sum_j a{a'}\0 \o {a'}\1$
(cf.\ right D2 quasibase eq.~(\ref{eq: right D2 quasibase})). 
This shows straightaway that the Galois mapping $\beta: A \o_B A
\to A \o_R T$ is  bijective.
Then $A \o_R T \o_R T \cong A \o_B A \o_B A$ via $\Phi := (\id_A \o \beta^{-1})(\beta^{-1} \o \id_T)$,
so coassociativity $(\id_A \o \cop)\rho_R = (\rho_R \o \id_T)\rho_R$ follows
from 
$$\Phi(( \id \o \cop_T) \circ \rho_R) = \sum_{j,k} \gamma_j(a) u^1_j \o_B \gamma_k(u^2_j)u^1_k \o_B u^2_k =
\sum_k 1 \o \gamma_k(a)u^1_k \o u^2_k = 1 \o 1 \o a
$$
$$= \sum_{j,k} \gamma_k(\gamma_j(a))u^1_k \o u^2_ku^1_j \o u^2_k = 
\Phi( (\rho_R \o \id)\rho_R(a)). $$ 
We note that $\rho_R$ is right $R$-linear, since
$$ \rho_R(ar) = \sum_j \gamma_j(ar) \o u_j \stackrel{\beta^{-1}}{\longmapsto} 1 \o_B ar = \beta^{-1}(\rho_R(a)r) $$
since $\rho_R(a)r = \sum_j \gamma_j(a) \o u_jr$. 
Also, $a\0 \eps_T(a\1) = $ $ \sum_j \gamma_j(a)u^1_j u^2_j  = a$
for all $a \in A$.  

Next, 
$$ \beta^{-1} (r \cdot a\0 \o a\1) = \sum_j r\gamma_j(a)  u_j = r \o_B a =
\sum_j \gamma_j(a)u^1_j r \o u^2_j = \beta^{-1}(a\0 \o t_R(r) a\1). $$
Whence the statement $\rho_R(aa') = \rho_R(a)\rho_R(a')$
makes sense for all $a,a' \in A$. 
We check the statement:  
$$\beta^{-1}(\rho_R(a)\rho_R(a')) = \sum_{j,k} \gamma_j(a) \gamma_k(a') u_j u_k =
 \sum_{j,k} \gamma_j(a) \gamma_k(a') u_k^1 u_j^1 \o u_j^2 u_k^2 
$$
$$ = 1 \o aa' = \sum_j \gamma_j(aa')u_j = \beta^{-1}(\rho_R(aa')). $$
Also $\rho_R(1_A) = 1_A \o_R 1_T$ since $\gamma_j(1_A) \in R$.
Finally we note that for each $b\in B$
$$ \rho_R(b) = \sum_j \gamma_j(b) \o_R u_j = b \o \sum_j \gamma_j(1)u_j = b \o 1_T $$
so $B \subseteq A^{\rm co \, \rho_R}$.  
Conversely, if $\rho_R(x) = x \o 1_T$ $ = \sum_j \gamma_j(x) \o u_j$ applying $\beta^{-1}$ we
obtain $x \o_B 1 = 1 \o_B x$. Let $f \in \End A_B$.   Then
applying $\mu (f \lambda(a) \o \id)$ to this we obtain 
 $f(ax) = f(a)x$  since $\lambda(a) \in \End A_B$
for each $a\in A$. It follows that $f \rho(x) = \rho(x) f$
so $\rho(x) \in \End {}_{\mathcal{E}}A$.  Since $A_B$ is balanced,
$\rho(x) = \rho(b)$ for some $ b\in B$, whence $x = b \in B$.  

The second part of the theorem is proven similarly. In the $\Leftarrow$ direction,
we convert the right $R$-bialgebroid $T$ to a left $R$-bialgebroid $T^{\rm op}$
with $s_L = t_R$, $t_L = s_L$, the same $R$-coring structure and opposite multiplication, which leads to  
the left $R$-$R$-bimodule structure coinciding with the usual $R$-$R$-bimodule structure
on $T$ in eq.~(\ref{eq: bimod}). We then define a left $T^{\rm op}$-comodule algebra structure on $A$
via $\rho_L: A \to T \o_R A$ defined via  left D2 quasibases by  
\begin{equation}
\rho_L(a) = a\-1 \o a\0 := \sum_i t_i \o \beta_i(a).
\end{equation}
The isomorphism $T \o_R A \stackrel{\cong}{\longrightarrow}$ $ A \o_B A$ given
by $t \o a \mapsto t^1 \o t^2a$ is inverse to the Galois mapping 
$\beta_L(a \o a') = $ $a\-1 \o a\0 a'$ by the left D2 quasibase eq.~(\ref{eq: left D2 quasibase}).
One needs the opposite multiplication of $T$ when showing $\rho_L(aa') = \rho_L(a) \rho_L(a')$
for $a,a'\in A$.  
\end{proof}
Let $T$ be a left finite projective right bialgebroid over some algebra $R$ in the next corollary.
 
\begin{cor}
Suppose $A \| B$ is a right $T$-extension.  If the Galois mapping $\beta$ is a split monic,
then $A \| B$ is a right $(A \o_B A)^B$-Galois extension.
\end{cor}
\begin{proof}
This follows from ${}_A A \o_B A_B \oplus * \cong {}_A A_B \o_R T$ and the arguments in the first few paragraphs
of the proof above (the balance argument makes only use of $A \| B$ being a right $T$-extension).  
Hence, $A \| B$ is right D2 and right balanced.  Whence $A \| B$ is a right Galois extension w.r.t.
the bialgebroid $(A \o_B A)^B$.
\end{proof}
Notice that $T$ is possibly not isomorphic to $(A \o_B A)^B$.  For example, one might start with 
a Hopf algebra Frobenius extension with split monic Galois map and conclude it is a weak Hopf-Galois extension
(if the centralizer is separable, the antipode being constructible from the Frobenius structure). 

We also observe that putting Theorems~\ref{th-characterization}
and \cite[2.6]{nov} together yields a type of endomorphism
ring theorem for depth two extensions, without a Frobenius
extension hypothesis (cf.\ \cite[Theorem 6.6]{KS}).

\begin{cor}
Suppose $A \| B$ is a depth two algebra extension.  Then
$\End {}_BA$ is a left D2 and left balanced extension of $A^{\rm op}$.
\end{cor}
\begin{proof}
In \cite[Theorem 2.6]{nov} it is established that $\End {}_BA$
is a left $S$-Galois extension of $\rho(A)$ $= \{ \rho(a) | a \in A \}$ where $\rho(a)(x) = xa$ for all $x,a \in A$.  But $S$
is a left and right f.g.\ projective left bialgebroid by the
D2 hypothesis (cf.\ proposition~\ref{prop-ess} below).  
It follows from the second statement in
the theorem above that $\End {}_BA \| A^{\rm op}$ is
left D2 and left balanced. 
\end{proof}

By carefully checking the proof of \cite[2.6]{nov}
for an airtight reliance on only right D2 quasibases, and referring to the proposition below, I believe it is likely that the corollary may be somewhat improved to show that a right D2 algebra extension has a left D2
left endomorphism algebra extension.

\section{Appendix}
In this section we answer some natural questions about the theory of one-sided depth two extensions.
One of the apparent questions after a reading of proposition~\ref{prop-tee} would
be if the  endomorphism algebra $S$ is also a bialgebroid over the centralizer, to which the next proposition 
provides an answer in the affirmative.
\begin{prop}
\label{prop-ess}
Suppose $A \| B$ is either a right or a left D2 extension with centralizer $R$.
Then $S$ is a left $R$-bialgebroid, which is either
right f.g.\ $R$-projective or left f.g.\ $R$-projective
respectively.
\end{prop}
\begin{proof}
The algebra structure comes from the usual composition of endomorphisms in $S =
\End {}_BA_B$.  The source and target mappings are $s_L(r) = \lambda(r)$
and $t_L(r) = \rho(r)$, whence the structure ${}_RS_R$ is given
by $$ r \cdot \alpha \cdot r' = \lambda(r) \rho(r') \alpha = r \alpha(-) r'.
$$

Suppose now we are given a right D2 structure on $A \| B$ by quasibases
$u_j \in T$, $\gamma_j \in S$.  
The $R$-coring structure on ${}_RS_R$ is given by a coproduct $\cop: S  \to S\o_R S$ defined by 
\begin{equation}
\label{eq: coproduct}
\cop(\alpha) = \sum_j \gamma_j \o u^1_j \alpha(u^2_j -) ,
\end{equation}
and a counit $\eps: S \to R$ given by
\begin{equation}
\label{eq: counit}
\eps(\alpha) = \alpha(1_A)
\end{equation}
Clearly $\eps$ is an $R$-$R$-bimodule mapping with $\eps(1_S) = 1_A$,
satisfying the counitality equations and $$\eps(\alpha \beta) = \eps(\alpha s_L(\eps(\beta))) = \eps(\alpha t_L(\eps(\beta))).$$ Also $\cop$ is right $R$-linear and $\cop(1_S) = 1_S \o_R 1_S$. 
By
making the identification 
$$ S \o_R S \cong \Hom ({}_B A \o_B A_B, {}_B A_B), \ \ \ \alpha \o \beta \longmapsto
(a \o a' \mapsto \alpha(a) \beta(a') )$$
with inverse $F \mapsto \sum_j \gamma_j \o u^1_j F(u^2_j \o - )$, 
we see that the coproduct is left $R$-linear, satisfies $\alpha\1 t_L(r) \o \alpha\2 =$ $ \alpha\1 \o \alpha\2 s_L(r)$
for all $r \in R$, and $\cop(\alpha \beta) = \cop(\alpha)\cop(\beta)$ for all $\alpha, \beta \in S$.  For with the independent variables $x,x' \in A$, $\alpha, \beta \in S$
and $r \in R$,
each of these expressions becomes equal  to $ r\alpha(xx')$, $\alpha(xrx')$, and $ \alpha(\beta(xx'))$
respectively. 

The coproduct $\cop$ is coassociative since 
$$ S \o_R S \o_R S \stackrel{\cong}{\longrightarrow} \Hom ({}_BA \o_B A \o_B A_B, {}_BA_B), \ \ \ \ 
\alpha \o \beta \o \gamma \longmapsto ( x \o y \o z \mapsto \alpha(x)\beta(y) \gamma(z) ) $$
with inverse given by
\begin{equation}
F \mapsto \sum_{i,j,k} \gamma_i \o u^1_i \gamma_j(u^2_i \gamma_k(-)) \o u^1_j F(u^2_j u^1_k \o u^2_k \o -) 
\end{equation}
Applying this identification to $(\cop \o \id_S) \cop(\alpha)$ and to $(\id_S \o \cop)\cop(\alpha)$
then to $x \o_B y \o_B z$ both expressions equal $\alpha(xyz)$. 

$S_R$ is f.g.\ projective since for each $\alpha \in S$, we have $\alpha = \sum_j \gamma_j h_j(\alpha)$
where $h_j \in \Hom (S_R, R_R)$ is defined by $h_j(\alpha) = u_j^1 \alpha(u_j^2)$. 

The proof that given left D2 quasibases $t_i \in T$, $\beta_i \in S$,
we have left f.g.\ projective left bialgebroid $S$ with identical bialgebroid structure is similar
and therefore omitted.
\end{proof}

Suppose $A \| B$ is right D2.  Then we have seen that $S$ is a right finite projective left bialgebroid
while $T$ is a left finite projective right bialgebroid.  There is a nondegenerate pairing between
$S$ and $T$ with values in the centralizer $R$ given by $\bra t \| \alpha \ket := $ $t^1 \alpha(t^2)$,
since \begin{equation}
\label{eq: pair}
\eta: {}_RT \stackrel{\cong}{\longrightarrow} \Hom (S_R, R_R) 
\end{equation}
via $\eta(t) = \bra t \| - \ket$ with inverse $\phi \mapsto$ $\sum_j \phi(\gamma_j)u_j$. 
By proposition \cite[2.5]{KS} a right f.g.\ projective left bialgebroid $S$ has a right $R$-bialgebroid
$R$-dual $S^*$.  The question is then if the  bialgebroid $S^*$ is isomorphic to the bialgebroid $T$
via $\eta$?  The question is partly answered in the affirmative by corollary \cite[5.3]{KS},
where it is shown without using left D2 quasibases that $T$ and $S^*$ are isomorphic via the pairing above
as algebras and $R$-$R$-bimodules.
\begin{cor}
Suppose $A \| B$ is right D2.  Then $T$ is isomorphic as right bialgebroids over $R$ to the right
$R$-dual of $S$ via $\eta$.  If $A \| B$ is left D2, then $T$ is isomorphic to
the bialgebroid left $R$-dual of $S$.
\end{cor}
\begin{proof}
What remains to check in the first statement is that $\eta$ is a homomorphism of $R$-corings
using right D2 quasibases.
We compute:
\begin{eqnarray*} 
\bra t\1 \cdot \bra t\2 \| \alpha' \ket \| \alpha \ket & = & \sum_j \bra t^1 \o_B \gamma_j(t^2)u^1_j \alpha'(u^2_j) \| \alpha \ket \\
& = & t^1 \alpha (\alpha'(t^2)) = \bra t \| \alpha \circ \alpha' \ket,
\end{eqnarray*}
Whence $\cop(\eta(t)) = \eta(t\1) \o \eta(t\2)$ by uniqueness \cite[2.5 (41)]{KS}. 

The proof of the last statement is similar to the first in using the pairing $[\alpha \| t ] :=$ $ \alpha(t^1)t^2$
and the right bialgebroid of the left dual of a left bialgebroid in \cite[2.6]{KS}. The details are left
to the reader.
\end{proof}



\begin{thebibliography}{XXXXXX}
\begin{small}
\bibitem{BaS}{I. B\'{a}lint and K.~Szlach\'anyi,
Finitary Galois extensions over non-commutative bases, KFKI preprint (2004), \texttt{RA/0412122}.}
\bibitem{Bo}{G.~B\"{o}hm,
Galois theory for Hopf algebroids, KFKI preprint (2004), \texttt{RA/0409513}.}
\bibitem{BW}{T.~Brzezi\'nski and R.~Wisbauer,
\textit{Corings and Comodules}, LMS \textbf{309}, Cambridge University Press, 2003.}
\bibitem{CDG}{S.~Caenepeel and E.~De Groot,
Galois theory for weak Hopf algebras, VUB preprint (2004), \texttt{RA/0406186}.}
\bibitem{Karl}{L.~Kadison,
Depth two and the Galois coring, preprint (2004), \texttt{RA/0408155}.}
\bibitem{AMA}{L.~Kadison,
An action-free characterization of weak Hopf-Galois extensions,
preprint (2004), \texttt{QA/0409589}.}
\bibitem{nov}{L.~Kadison,
Normal Hopf subalgebras, depth two and Galois extensions, preprint (2004), 
\texttt{QA/0411129}.} 
\bibitem{KN}{L.~Kadison and D.~Nikshych,
Hopf algebra actions on strongly separable extensions of depth two,
\textit{Adv.\ in Math.} \textbf{163} (2001), 258--286.}
\bibitem{KS}{L.~Kadison and K.~Szlach\'anyi,
Bialgebroid actions on depth two extensions and duality, \textit{Adv.\
in Math.} \textbf{179} (2003), 75--121. 
\texttt{RA/0108067}}
\bibitem{Mo}{S.~Montgomery,
``Hopf Algebras and Their Actions on Rings,''  CBMS Regional Conf.\ 
Series in Math.\ 
Vol.\ 82, AMS, Providence, 1993.}
\end{small}
\end{thebibliography}
\end{document}